# THE SECOND HILBERT COEFFICIENTS AND THE SECTIONAL GENERA OF IDEALS

KAZUHO OZEKI

## 1. Introduction

The purpose of our paper is to study the relationship between the second Hilbert coefficients and the sectional genera of ideals.

Let $A$ be a Cohen-Macaulay local ring with maximal ideal $\mathfrak{m}$ and $d = \dim A > 0$. For simplicity, throughout this paper, we assume that the residue class field $A/\mathfrak{m}$ of $A$ is infinite. Let $I$ be a fixed $\mathfrak{m}$-primary ideal in $A$ and $Q = (a_1, a_2, \cdots, a_d)$ a parameter ideal of $A$ which forms a reduction of $I$, that is the equality $I^{n+1} = QI^n$ holds true for some $n > 0$. Let $\ell_A(N)$ denote, for an $A$-module $N$, the length of $N$. Then there exist integers $\{\mathrm{e}_i(I)\}_{0 \le i \le d}$ such that

$$\ell_A(A/I^{n+1}) = \mathrm{e}_0(I)\binom{n+d}{d} - \mathrm{e}_1(I)\binom{n+d-1}{d-1} + \cdots + (-1)^d \mathrm{e}_d(I)$$

for all $n \gg 0$. We call $\mathrm{e}_i(I)$ the $i$-th Hilbert coefficient of $I$ and especially call the leading coefficient $\mathrm{e}_0(I) = \mathrm{e}_0(Q) = \ell_A(A/Q)$ $(> 0)$ the multiplicity of $I$. Let

$$\mathrm{R}'(I) = A[It, t^{-1}] \subseteq A[t, t^{-1}] \quad \text{and} \quad \mathrm{gr}(I) = \mathrm{R}'(I)/t^{-1}\mathrm{R}'(I)$$

where $t$ is an indeterminate over $A$.

As a classical result of Northcott [No], the inequality $\mathrm{e}_1(I) \ge \mathrm{e}_0(I) - \ell_A(A/I)$ $(\ge 0)$ holds true for every $\mathfrak{m}$-primary ideals $I$ in $A$, and Huneke [H] and Ooishi [O] showed that the equality $\mathrm{e}_1(I) = \mathrm{e}_0(I) - \ell_A(A/I)$ holds true if and only if $I^2 = QI$. When this is the case, the associated graded ring $\mathrm{gr}(I)$ and the fiber cone $\mathrm{F}(I) = \bigoplus_{n \ge 0} I^n/\mathfrak{m}I^n$ of $I$ are both Cohen-Macaulay, and the Rees algebra $\mathrm{R}(I) = A[It] \subseteq A[t]$ of $I$ is also a Cohen-Macaulay ring, provided $d \ge 2$. We also notice that, Kirby and Mehran [KM] were able to show that $\mathrm{e}_1(I) \le \binom{\mathrm{e}_0(I)}{2}$.

The purpose of this paper is to study the second Hilbert coefficients $\mathrm{e}_2(I)$ of $I$. As is well known, Narita [Na] showed that the inequality $\mathrm{e}_2(I) \ge 0$ holds true for any $\mathfrak{m}$-primary ideal $I$ in $A$. In [KM], it was proved that for $\mathfrak{m}$-primary ideals $I$ in $A$, an upper bound

$$\mathrm{e}_2(I) \le \binom{\mathrm{e}_1(I) + 1}{2}$$

of the second Hilbert coefficient $\mathrm{e}_2(I)$ in terms of the first Hilbert coefficient $\mathrm{e}_1(I)$.

Our first purpose is to give an upper bound of the second Hilbert coefficients $\mathrm{e}_2(I)$ in terms of the sectional genera of $I$. Let

$$\mathrm{g}_s(I) = \ell_A(A/I) - \mathrm{e}_0(I) + \mathrm{e}_1(I)$$

denotes the sectional genus of $I$. When $d = \dim A \geq 2$, we shall show that the inequality

$$\mathrm{e}_2(I) \leq \binom{\mathrm{g}_s(I)+1}{2}$$

holds true for every $\mathfrak{m}$-primary ideal $I$ in $A$ (Proposition 3.2). The following theorem which is the first main result of this paper shows that the upper bound $\mathrm{e}_2(I) \leq \binom{\mathrm{g}_s(I)+1}{2}$ is sharp, clarifying when the equality $\mathrm{e}_2(I) = \binom{\mathrm{g}_s(I)+1}{2}$ holds true.

**Theorem 1.1.** *Suppose that $d \geq 2$ and let $g = \mathrm{g}_s(I)$. Then the following three conditions are equivalent:*

(1) $\mathrm{e}_2(I) = \binom{g+1}{2}$,
(2) $\ell_A(I^2/QI) \leq 1$,
(3) $\ell_A(I^{k+1}/QI^k) = 1$ *for all* $1 \leq k \leq g$ *and* $I^{g+2} = QI^{g+1}$.

*When this is the case we have the following.*

(i) $\mathrm{e}_i(I) = \binom{g+1}{i}$ *for* $3 \leq i \leq d$,
(ii) $\operatorname{depth} \mathrm{gr}(I) \geq d-1$, *and*
(iii) $\mathrm{gr}(I)$ *is a Cohen-Macaulay ring if and only if* $Q \cap I^2 = QI$ *and* $I^3 = QI^2$.

Thus the second Hilbert coefficients $\mathrm{e}_2(I)$ bounded above by $\binom{\mathrm{g}_s(I)+1}{2}$. It seems now natural to ask what happens on the $\mathfrak{m}$-primary ideals $I$ with $\mathrm{e}_2(I) < \binom{\mathrm{g}_s(I)+1}{2}$. The second main result of this paper answers the question and is stated as follows (Theorem 3.4).

**Theorem 1.2.** *Suppose that $d \geq 2$. Let $g = \mathrm{g}_s(I)$ and assume that $g \geq 2$. Then the following three conditions are equivalent:*

(1) $\mathrm{e}_2(I) = \binom{g}{2} + 1$,
(2) $\binom{g}{2} + 1 \leq \mathrm{e}_2(I) < \binom{g+1}{2}$,
(3) $\ell_A(I^2/QI) = 2$, $\ell_A(I^{k+1}/QI^k) = 1$ *for all* $2 \leq k \leq g-1$, *and* $I^{g+1} = QI^g$.

*When this is the case, we have the following*:

(i) $\mathrm{e}_i(I) = \binom{g}{i}$ *for* $3 \leq i \leq d$,
(ii) $\operatorname{depth} \mathrm{gr}(I) \geq d-1$, *and*
(iii) $\mathrm{gr}(I)$ *is a Cohen-Macaulay ring if and only if* $Q \cap I^2 = QI$, $Q \cap I^3 = QI^2$, *and* $I^4 = QI^3$.

As a direct consequence of Theorem 1.2, we get $\mathrm{e}_2(I) = \binom{\mathrm{g}_s(I)}{2} + 1$ or $\mathrm{e}_2(I) = \binom{\mathrm{g}_s(I)+1}{2}$, if $\mathrm{e}_2(I) \geq \binom{\mathrm{g}_s(I)}{2} + 1$ (Corollary 4.2).

We now briefly explain how this paper is organized. In Section 2 we will summarize, for the later use in this paper, some auxiliary results on the Hilbert coefficients and the sectional genera of $\mathfrak{m}$-primary ideals. We shall prove Theorem 1.1 in Section 3 (Theorem 3.4). Theorem 1.2 will be proven in Section 4 (Theorem 4.1). We will show in Section 5 examples of $\mathfrak{m}$-primary ideals satisfying the equality in Theorem 1.1 (1) and Theorem 1.2 (1), respectively.

In what follows, unless otherwise specified, let $A$ be a Cohen-Macaulay local ring with maximal ideal $\mathfrak{m}$ and $d = \dim A > 0$. We throughout assume that the field $A/\mathfrak{m}$ is infinite. Let $I$ be an $\mathfrak{m}$-primary ideal in $A$ and $Q = (a_1, a_2, \cdots, a_d)$ a parameter ideal of $A$ which forms a reduction of $I$. Let $\mathrm{r}_Q(I) = \sup\{n \geq 0 \mid I^{n+1} = QI^n\}$ denotes the reduction number of $I$ with respect to $Q$. For each $\mathfrak{m}$-primary ideal $I$ in $A$ we set

$$R' = \mathrm{R}'(I) = A[It, t^{-1}] \subseteq A[t, t^{-1}] \ \text{ and } \ \mathrm{gr}(I) = R'/t^{-1}R'$$

where $t$ is an indeterminate over $A$.

## 2. Preliminary steps

In this section we summarize some basic properties and known results of the Hilbert coefficients, and the sectional genera of ideals, which we need throughout this paper.

The following result is, more or less, known (c.f. [CPP, Corollary 1.3], [RV2, Section 2.2]). Let us indicate a brief proof for the sake of completeness, because it plays a key role in the proofs of our main theorems. Set $\mathrm{H}(\mathrm{gr}(I), t)$ denotes the Hilbert series of the associated graded ring $\mathrm{gr}(I)$ of $I$.

**Lemma 2.1.** *Suppose that $d = 1$. Then we have*

$$\mathrm{H}(\mathrm{gr}(I), t) = \frac{\ell_A(A/I) + \sum_{k=1}^{r-1}\{\ell_A(I^k/QI^{k-1}) - \ell_A(I^{k+1}/QI^k)\}t^k + \ell_A(I^r/QI^{r-1})t^r}{1-t}$$

*where $r = \mathrm{r}_Q(I)$.*

*Proof.* Let

$$\mathrm{H}(\mathrm{gr}(I), t) = \frac{\sum_{k \geq 0} h_k t^k}{1-t}$$

with $h_k \in \mathbb{Z}$ for $k \geq 0$. Then we have $h_0 = \ell_A([\mathrm{gr}(I)]_0) = \ell_A(A/I)$ and

$$h_k = \ell_A([\mathrm{gr}(I)]_k) - \ell_A([\mathrm{gr}(I)]_{k-1}) = \ell_A(I^k/I^{k+1}) - \ell_A(I^{k-1}/I^k)$$

for all $k \geq 1$. Thus, by the exact sequences

$$0 \to I^{k+1}/QI^k \to I^k/QI^k \to I^k/I^{k+1} \to 0$$

and

$$0 \to I^{k-1}/I^k \xrightarrow{a_1} I^k/QI^k \to I^k/QI^{k-1} \to 0$$

of $A$-modules, we get

$$\begin{aligned} h_k &= \ell_A(I^k/I^{k+1}) - \ell_A(I^{k-1}/I^k) \\ &= \{\ell_A(I^k/QI^k) - \ell_A(I^{k+1}/QI^k)\} - \{\ell_A(I^k/QI^k) - \ell_A(I^k/QI^{k-1})\} \\ &= \ell_A(I^k/QI^{k-1}) - \ell_A(I^{k+1}/QI^k) \end{aligned}$$

for all $k \geq 1$ as required. □

We also need the following lemma.

**Lemma 2.2.** *Suppose that $d \geq 2$ and let $a \in I$ be a superficial element for $I$. Then we have $\mathrm{g}_s(I) = \mathrm{g}_s(I/(a))$.*

*Proof.* We set $\overline{A} = A/(a)$ and $\overline{I} = I/(a)$. Then we have $\mathrm{g}_s(\overline{I}) = \ell_A(\overline{A}/\overline{I}) - \mathrm{e}_0(\overline{I}) + \mathrm{e}_1(\overline{I}) = \ell_A(A/I) - \mathrm{e}_0(I) + \mathrm{e}_1(I) = g_s(I)$ as $\mathrm{e}_i(\overline{I}) = \mathrm{e}_i(I)$ for all $i = 0, 1$, and $\ell_A(\overline{A}/\overline{I}) = \ell_A(A/I)$. □

Thanks to Lemma 2.1, the following result holds true.

**Lemma 2.3.** *Suppose that $d = 2$ and let $a \in I$ be a superficial element for $I$. Then we have the following, where $\overline{A} = A/(a)$, $\overline{I} = I/(a)$, $\overline{Q} = Q/(a)$, $r' = \mathrm{r}_{\overline{Q}}(\overline{I})$, and $m$ be a positive integer such that $I^{n+1} :_A a = I^n$ for all $n \geq m$.*

(1) $\mathrm{g}_s(I) = \displaystyle\sum_{k=1}^{r'-1} \ell_A(\overline{I}^{k+1}/\overline{Q}\,\overline{I}^k)$,

(2) $\mathrm{e}_2(I) = \displaystyle\sum_{k=1}^{r'-1} k{\cdot}\ell_A(\overline{I}^{k+1}/\overline{Q}\,\overline{I}^k) - \sum_{k=1}^{m-1} \ell_A([I^{k+1} :_A a]/I^k)$.

*Proof.* We set $v_k = \ell_A(\overline{I}^{k+1}/\overline{Q}\,\overline{I}^k)$ for all $k \geq 0$. Then, by Lemma 2.1, we have

$$\mathrm{H}(\mathrm{gr}(\overline{I}), t) = \frac{\ell_A(\overline{A}/\overline{I}) + \sum_{k=1}^{r'-1}\{v_{k-1} - v_k\}t^k + v_{r'-1}t^{r'}}{1-t}$$

so that $\mathrm{e}_1(\overline{I}) = \sum_{k=0}^{r'-1} v_k$ and $\mathrm{e}_2(\overline{I}) = \sum_{k=1}^{r'-1} k{\cdot}v_k$. Therefore, we get

$$\mathrm{g}_s(I) = \mathrm{g}_s(\overline{I}) = \ell_A(\overline{A}/\overline{I}) - \mathrm{e}_0(\overline{I}) + \mathrm{e}_1(\overline{I}) = \sum_{k=1}^{r'-1} v_k$$

by Lemma 2.2, because $v_0 = \ell_A(\overline{I}/\overline{Q}) = \mathrm{e}_0(\overline{I}) - \ell_A(\overline{A}/\overline{I})$. We also get

$$\mathrm{e}_2(I) = \mathrm{e}_2(\overline{I}) - \sum_{k=1}^{m-1} \ell_A([I^{k+1} :_A a]/I^k) = \sum_{k=1}^{r'-1} k{\cdot}v_k - \sum_{k=1}^{m-1} \ell_A([I^{k+1} :_A a]/I^k)$$

as required (c.f. [RV2, Proposition 1.2]). □

The following result is due to Huckaba and Marley [HM].

**Theorem 2.4.** ([HM, Theorem 4.7]) *Suppose that* $d \geq 1$*. Then we have the inequality*

$$\mathrm{e}_1(I) \leq \sum_{k \geq 0} \ell_A(I^{k+1}/QI^k)$$

*and the following two conditions are equivalent:*

(1) $\mathrm{e}_1(I) = \sum_{k \geq 0} \ell_A(I^{k+1}/QI^k)$,

(2) $\operatorname{depth} \operatorname{gr}(I) \geq d - 1$.

*When this is the case we have*

$$\mathrm{e}_i(I) = \sum_{k \geq i-1} \binom{k}{i-1} \ell_A(I^{k+1}/QI^k)$$

*for all* $2 \leq i \leq d$.

*Proof.* See, for example, [RV2, Theorem 2.5]. $\square$

The following result was proved by Wang [W] and, at the same time by Rossi and Valla [RV1].

**Theorem 2.5.** ([RV1, W]) *Suppose that* $d \geq 1$ *and assume that* $\ell_A(I^2/QI) \leq 1$*. Then we have the following.*

(1) $\operatorname{depth} \operatorname{gr}(I) \geq d - 1$,

(2) $\mathrm{e}_i(I) = \binom{\mathrm{g}_s(I)+1}{i}$ *for all* $2 \leq i \leq d$.

*Proof.* (1) See [S1, W] ([RV2, Theorem 4.4] also).

(2) We have $\ell_A(I^{k+1}/QI^k) \leq 1$ for all $k \geq 1$ by [CPP, Lemma 3.8] (see [S1, W] also). Since $\operatorname{depth} \operatorname{gr}(I) \geq d-1$, we have $\mathrm{g}_s(I) = \ell_A(A/I) - \mathrm{e}_0(I) + \mathrm{e}_1(I) = \sum_{k \geq 1} \ell_A(I^{k+1}/QI^k)$ by Theorem 2.4, so that we have $\ell_A(I^{k+1}/QI^k) = 1$ for all $1 \leq k \leq g$ and $I^{g+2} = QI^{g+1}$. Thus, we get $\mathrm{e}_i(I) = \binom{g+1}{i}$ for all $2 \leq i \leq d$ by Theorem 2.4. $\square$

Thus the proofs of the implication (2) $\Rightarrow$ (1) and the last assertion ($i$) in Theorem 1.1 are given by Theorem 2.5.

## 3. An upper bound for the second Hilbert coefficient

The purpose of this section is to estimate the second Hilbert coefficients of $\mathfrak{m}$-primary ideals in terms of the sectional genera.

Let us begin with the following lemma.

**Lemma 3.1.** *Let* $\ell \geq 0$ *be an integer. Suppose that* $\{v_k\}_{k \geq 1}$ *is the set of integers such that*

($i$) $v_k \geq 0$ *for all* $k \geq 1$,

($ii$) $\sum_{k \geq 1} v_k \leq \ell$, *and*

($iii$) $v_j = 0$ *for all* $j \geq k$ *once* $v_k = 0$ *for some* $k \geq 1$.

*Then we have*

$$\sum_{k\geq 1} k{\cdot}v_k \leq \binom{\ell+1}{2}$$

*and the following two conditions are equivalent:*

(1) $\sum_{k\geq 1} k{\cdot}v_k = \binom{\ell+1}{2}$,
(2) $v_k = 1$ *for all* $1 \leq k \leq \ell$, *and* $v_k = 0$ *for all* $k \geq \ell+1$.

*Proof.* We proceed by induction on $\ell$. When $\ell \leq 1$, we have nothing to do. Assume that $\ell \geq 2$ and that our assertion holds true for $\ell - 1$. We may assume that $v_1 \geq 1$. We set $w_k = v_{k+1}$ for all $k \geq 1$. Then, by the hypothesis of induction on $\ell$, we have

$$\sum_{k\geq 2}(k-1){\cdot}v_k = \sum_{k\geq 1} k{\cdot}w_k \leq \binom{\ell - v_1 + 1}{2} \leq \binom{\ell}{2}$$

because $\sum_{k\geq 1} w_k = \sum_{k\geq 1} v_k - v_1 \leq \ell - v_1 \leq \ell - 1$. Therefore we get the required inequality

$$\sum_{k\geq 1} k{\cdot}v_k = \sum_{k\geq 2}(k-1){\cdot}v_k + \sum_{k\geq 1} v_k \leq \binom{\ell}{2} + \ell = \binom{\ell+1}{2}.$$

In the rest of our proof of Lemma 3.1, we have only to show the implication (1) $\Rightarrow$ (2). Assume that $\sum_{k\geq 1} k{\cdot}a_k = \binom{\ell+1}{2}$. Then we have

$$\sum_{k\geq 1} k{\cdot}w_k = \binom{\ell - v_1 + 1}{2} = \binom{\ell}{2}$$

by the above argument, and hence $v_1 = 1$. Thus, by the hypothesis of induction on $\ell$, we get $v_k = w_{k-1} = 0$ for all $2 \leq k \leq \ell$ and $v_k = w_{k-1} = 0$ for all $k \geq \ell + 1$ as required. $\square$

We notice that the set of integers $\{\ell_A(I^{k+1}/QI^k)\}_{k\geq 1}$, for an $\mathfrak{m}$-primary ideal $I$ and a reduction $Q$ of $I$, satisfies all conditions $(i)$, $(ii)$, and $(iii)$ in Lemma 3.1.

The following result gives an upper bound for $\mathrm{e}_2(I)$ in terms of the sectional genus $\mathrm{g}_s(I)$ of $I$, and also shows that the implication (1) $\Rightarrow$ (3) and the assertion $(ii)$ in Theorem 1.1 (Theorem 3.4) hold true.

**Proposition 3.2.** *Suppose that* $d \geq 2$ *and let* $g = \mathrm{g}_s(I)$. *Then we have the inequality*

$$\mathrm{e}_2(I) \leq \binom{g+1}{2}.$$

*If* $\mathrm{e}_2(I) = \binom{g+1}{2}$ *then the following two assertions hold true.*

(1) depth $\mathrm{gr}(I) \geq d - 1$,
(2) $\ell_A(I^{k+1}/QI^k) = 1$ *for all* $1 \leq k \leq g$ *and* $I^{g+2} = QI^{g+1}$.

Thanks to Proposition 3.2, we get the following inequality which is given by Kirby and Mehran [KM].

**Corollary 3.3.** *Suppose that $d \geq 2$. Then we have the inequality*

$$\mathrm{e}_2(I) \leq \binom{\mathrm{e}_1(I)+1}{2}.$$

*We furthermore have $I = Q$, once the equality $\mathrm{e}_2(I) = \binom{\mathrm{e}_1(I)+1}{2}$ holds true.*

*Proof.* The inequality $\mathrm{e}_2(I) \leq \binom{\mathrm{e}_1(I)+1}{2}$ holds true by Proposition 3.2. Assume that $\mathrm{e}_2(I) = \binom{\mathrm{e}_1(I)+1}{2}$. Then, because $\mathrm{e}_2(I) = \binom{\mathrm{e}_1(I)+1}{2} = \binom{\mathrm{g}_s(I)+1}{2}$, we have $\mathrm{e}_1(I) = \mathrm{g}_s(I)$, so that $\ell_A(I/Q) = \mathrm{e}_0(I) - \ell_A(A/I) = \mathrm{e}_1(I) - \mathrm{g}_s(I) = 0$. Thus we get $I = Q$ as required. $\square$

*Proof of Proposition 3.2.* We proceed by induction on $d$. Since the residue class field $A/\mathfrak{m}$ of $A$ is infinite, we may choose an element $a \in Q \backslash \mathfrak{m}Q$ is superficial for $I$. We set $\overline{A} = A/(a)$, $\overline{I} = I/(a)$, $\overline{Q} = Q/(a)$, and $r' = \mathrm{r}_{\overline{Q}}(\overline{I})$.

Suppose that $d = 2$. Let $v_k = \ell_A(\overline{I}^{k+1}/\overline{Q}\,\overline{I}^k)$ for all $k \geq 1$. Then we have $g = \sum_{k=1}^{r'-1} v_k$ by Lemma 2.3. Therefore, by Lemma 2.3 and Lemma 3.1, we get

$$\mathrm{e}_2(I) = \sum_{k=1}^{r'-1} k{\cdot}v_k - \sum_{k=1}^{m-1} \ell_A([I^{k+1} :_A a]/I^k) \leq \sum_{k=1}^{r'-1} k{\cdot}v_k \leq \binom{g+1}{2},$$

where $m \geq 1$ be an integer such that $I^{k+1} :_A a = I^k$ for all $k \geq m$. Assume $\mathrm{e}_2(I) = \binom{g+1}{2}$ then we have $\sum_{k=1}^{r'-1} k{\cdot}v_k = \binom{g+1}{2}$. Therefore, we get $\ell_A(\overline{I}^{k+1}/\overline{Q}\,\overline{I}^k) = v_k = 1$ for all $1 \leq k \leq g$ and $\ell_A(\overline{I}^{k+1}/\overline{Q}\,\overline{I}^k) = v_k = 0$ for all $k \geq g+1$ by Lemma 3.1, and hence we have $\overline{I}^{g+2} = \overline{Q}\,\overline{I}^{g+1}$. We also have $\operatorname{depth} \mathrm{gr}(I) \geq 1$ because $I^{k+1} :_A a = I^k$ for all $k \geq 1$. Then, since $at$ forms a $\mathrm{gr}(I)$-regular element, we get $\ell_A(I^{k+1}/QI^k) = 1$ for all $1 \leq k \leq g$ and $I^{g+2} = QI^{g+1}$ as required.

Assume that $d \geq 3$ and that our assertion holds true for $d-1$. Then by the hypothesis of induction on $d$, we get the required inequality

$$\mathrm{e}_2(I) = \mathrm{e}_2(\overline{I}) \leq \binom{\mathrm{g}_s(\overline{I})+1}{2} = \binom{g+1}{2}$$

because $\mathrm{g}_s(\overline{I}) = g$ by Lemma 2.2. Assume that $\mathrm{e}_2(I) = \binom{g+1}{2}$ then we have $\mathrm{e}_2(\overline{I}) = \binom{\mathrm{g}_s(\overline{I})+1}{2}$. The hypothesis of induction on $d$ says that $\operatorname{depth} \mathrm{gr}(\overline{I}) \geq (d-1)-1 = d-2 > 0$, $\ell_A(\overline{I}^{k+1}/\overline{Q}\,\overline{I}^k) = 1$ for all $1 \leq k \leq g$, and $\overline{I}^{g+2} = \overline{Q}\,\overline{I}^{g+1}$. Then, thanks to Sally's technique ([S1], [HM, Lemma 2.2]), $at$ forms a $\mathrm{gr}(I)$-regular element. Thus we get $\operatorname{depth} \mathrm{gr}(I) \geq d-1$, $\ell_A(I^{k+1}/QI^k) = 1$ for all $1 \leq k \leq g$, and $I^{g+2} = QI^{g+1}$ as required. This completes the proof of Proposition 3.2. $\square$

Thus the second Hilbert coefficients $\mathrm{e}_2(I)$ of $\mathfrak{m}$-primary ideals $I$ are bounded above by $\binom{\mathrm{g}_s(I)+1}{2}$. It seems now natural to ask what happens on the $\mathfrak{m}$-primary ideals $I$ of $A$, once the equality $\mathrm{e}_2(I) = \binom{\mathrm{g}_s(I)+1}{2}$ is attained. The main result of this section answers the question and is stated as follows (Theorem 1.1).

**Theorem 3.4.** *Suppose that $d \geq 2$ and let $g = \mathrm{g}_s(I)$. Then the following three conditions are equivalent:*

(1) $\mathrm{e}_2(I) = \binom{g+1}{2}$,

(2) $\ell_A(I^2/QI) \leq 1$,

(3) $\ell_A(I^{k+1}/QI^k) = 1$ *for all* $1 \leq k \leq g$ *and* $I^{g+2} = QI^{g+1}$.

*When this is the case we have the following.*

(i) $\mathrm{e}_i(I) = \binom{g+1}{i}$ *for* $3 \leq i \leq d$,

(ii) $\operatorname{depth} \operatorname{gr}(I) \geq d-1$, *and*

(iii) $\operatorname{gr}(I)$ *is a Cohen-Macaulay ring if and only if* $Q \cap I^2 = QI$ *and* $I^3 = QI^2$.

*Proof.* (1) $\Rightarrow$ (3) and (*ii*) By Proposition 3.2.

(3) $\Rightarrow$ (2) It is clear.

(2) $\Rightarrow$ (1) and (*i*) By Theorem 2.5.

Now we have only to show that the assertion (*iii*) holds true. Assume that $\operatorname{gr}(I)$ is a Cohen-Macaulay ring. Then we have $Q \cap I^{n+1} = QI^n$ for all $n \geq 1$ by Valabrega-Valla criterion ([VV]). Since $\ell_A(I^2/QI) \leq 1$, we have $I^3 \subseteq QI$, so that $I^3 = Q \cap I^3 = QI^2$. The converse also holds true by Valabrega-Valla criterion ([VV]). This completes the proof of Theorem 3.4. □

From now on, we introduce two consequences of Theorem 3.4.

The following result gives a characterization of ideals $I$ with $\ell_A(I/Q) \leq 1$ in terms of $\mathrm{e}_1(I)$ and $\mathrm{e}_2(I)$.

**Corollary 3.5.** *Suppose that $d \geq 2$. Then the following two conditions are equivalent:*

(1) $\mathrm{e}_2(I) = \binom{\mathrm{e}_1(I)}{2}$,

(2) $\mathrm{e}_0(I) \leq \ell_A(A/I) + 1$, *that is* $\ell_A(I/Q) \leq 1$.

*When this is the case, we have the following.*

(*i*) $\mathrm{e}_i(I) = \binom{\mathrm{e}_1(I)}{i}$ *for* $3 \leq i \leq d$,

(*ii*) $\operatorname{depth} \operatorname{gr}(I) \geq d-1$, *and*

(*iii*) $\operatorname{gr}(I)$ *is a Cohen-Macaulay ring if and only if* $I^2 = QI$.

*Proof.* (1) $\Rightarrow$ (2) We may assume that $\ell_A(I/Q) \geq 1$, hence we have $\mathrm{g}_s(I) = \ell_A(A/I) - \mathrm{e}_0(I) + \mathrm{e}_1(I) = \mathrm{e}_1(I) - \ell_A(I/Q) \leq \mathrm{e}_1(I) - 1$. Therefore the inequalities

$$\mathrm{e}_2(I) \leq \binom{\mathrm{g}_s(I)+1}{2} \leq \binom{\mathrm{e}_1(I)}{2} = \mathrm{e}_2(I)$$

follow by Proposition 3.2. Then, since $\mathrm{e}_1(I) = \mathrm{g}_s(I) + 1$, we get $\ell_A(I/Q) = \mathrm{e}_0(I) - \ell_A(A/I) = \mathrm{e}_1(I) - \mathrm{g}_s(I) = 1$.

(2) $\Rightarrow$ (1) We may assume that $\ell_A(I/Q) = 1$. Let $I = Q + (x)$ with $x \in I$. Then $I^2 = QI + (x^2)$, so that $\ell_A(I^2/QI) \leq 1$ because $\mathfrak{m}I \subseteq Q$. Therefore, because $\mathrm{g}_s(I) =$

$\ell_A(A/I) - \mathrm{e}_0(I) + \mathrm{e}_1(I) = \mathrm{e}_1(I) - \ell_A(I/Q) = \mathrm{e}_1(I) - 1$, the equality $\mathrm{e}_2(I) = \binom{\mathrm{g}_s(I)+1}{2} = \binom{\mathrm{e}_1(I)}{2}$ holds true by Theorem 2.5.

(*i*) and (*ii*) By Theorem 1.1.

(*iii*) Assume that $\mathrm{gr}(I)$ is a Cohen-Macaulay ring. Then, because $\ell_A(I/Q) \leq 1$, we have $I^2 \subseteq Q$, so that $I^2 = Q \cap I^2 = QI$ by Theorem 3.4 (iii). The converse also follows by Theorem 3.4 (iii). □

In the end of this section we review a result of Sally [S1] in order to see how our Theorem 3.4 works to prove or improve them.

**Corollary 3.6.** *Suppose that $d \geq 2$. Assume that $\mathrm{g}_s(I) = 1$ and $\mathrm{e}_2(I) \neq 0$. Then the following assertions hold true.*

(1) $\mathrm{e}_2(I) = 1$ *and* $\mathrm{e}_i(I) = 0$ *for all* $3 \leq i \leq d$*,*
(2) $\ell_A(I^2/QI) = 1$ *and* $I^3 = QI^2$*, and*
(3) $\operatorname{depth} \mathrm{gr}(I) \geq d - 1$.

## 4. The second border

The second Hilbert coefficients $\mathrm{e}_2(I)$ of $\mathfrak{m}$-primary ideals $I$ are bounded above by $\binom{\mathrm{g}_s(I)+1}{2}$. In this section we explore the second border of $\mathrm{e}_2(I)$ in terms of the sectional genus $\mathrm{g}_s(I)$. The main result of this section is stated as follows.

**Theorem 4.1.** *Suppose that $d \geq 2$. Let $g = \mathrm{g}_s(I, M)$ and assume that $g \geq 2$. Then the following three conditions are equivalent:*

(1) $\mathrm{e}_2(I) = \binom{g}{2} + 1$,
(2) $\binom{g}{2} + 1 \leq \mathrm{e}_2(I) < \binom{g+1}{2}$,
(3) $\ell_A(I^2/QI) = 2$, $\ell_A(I^{k+1}/QI^k) = 1$ *for all* $2 \leq k \leq g-1$*, and* $I^{g+1} = QI^g$.

*When this is the case, we have the following.*

(i) $\mathrm{e}_i(I) = \binom{g}{i}$ *for* $3 \leq i \leq d$*,*
(ii) $\operatorname{depth} \mathrm{gr}(I) \geq d - 1$*, and*
(iii) $\mathrm{gr}(I)$ *is a Cohen-Macaulay ring if and only if* $Q \cap I^2 = QI$, $Q \cap I^3 = QI^2$*, and* $I^4 = QI^3$.

As a direct consequence of Theorem 4.1, we have the following.

**Corollary 4.2.** *Suppose that $d \geq 2$ and assume that $\mathrm{e}_2(I) \geq \binom{\mathrm{g}_s(I)}{2} + 1$. Then we have $\mathrm{e}_2(I) = \binom{\mathrm{g}_s(I)}{2} + 1$ or $\mathrm{e}_2(I) = \binom{\mathrm{g}_s(I)+1}{2}$, and $\operatorname{depth} \mathrm{gr}(I) \geq d - 1$.*

Before giving a proof of Theorem 4.1, let us begin with the following. We set

$$\widetilde{I} = \bigcup_{n \geq 0} [I^{n+1} :_A I^n]$$

denotes the Ratliff-Rush closure of $I$, which is the largest $\mathfrak{m}$-primary ideal in $A$ such that $I \subseteq \widetilde{I}$ and $\mathrm{e}_i(\widetilde{I}) = \mathrm{e}_i(I)$ for all $0 \leq i \leq d$ (c.f. [RR]).

**Lemma 4.3.** *Suppose that $d \geq 2$. We have $\widetilde{I} = I$, if $\mathrm{e}_2(I) \geq \binom{\mathrm{g}_s(I)}{2} + 1$.*

*Proof.* Assume that $\widetilde{I} \supsetneq I$. Then

$$\begin{aligned} \mathrm{g}_s(\widetilde{I}) &= \ell_A(A/\widetilde{I}) - \mathrm{e}_0(\widetilde{I}) + \mathrm{e}_1(\widetilde{I}) \\ &= \{\ell_A(A/I) - \ell_A(\widetilde{I}/I)\} - \mathrm{e}_0(I) + \mathrm{e}_1(I) \\ &= \mathrm{g}_s(I) - \ell_A(\widetilde{I}/I) \leq \mathrm{g}_s(I) - 1 \end{aligned}$$

because $\mathrm{e}_i(\widetilde{I}) = \mathrm{e}_i(I)$ for all $i = 0, 1$. Therefore we get

$$\mathrm{e}_2(I) = \mathrm{e}_2(\widetilde{I}) \leq \binom{\mathrm{g}_s(\widetilde{I}) + 1}{2} = \binom{\mathrm{g}_s(I) - \ell_A(\widetilde{I}/I) + 1}{2} \leq \binom{\mathrm{g}_s(I)}{2}$$

because $\mathrm{e}_2(\widetilde{I}) \leq \binom{\mathrm{g}_s(\widetilde{I})+1}{2}$ by Proposition 3.2. However, since $\mathrm{e}_2(I) \geq \binom{\mathrm{g}_s(I)}{2} + 1$ by our assumption, it is impossible, so that we can get a required contradiction. Thus $\widetilde{I} = I$ as required. $\square$

We are now in a position to prove Theorem 4.1.

*Proof of Theorem 4.1.* (1) $\Rightarrow$ (2) It is clear.

(2) $\Rightarrow$ (3) and (*ii*). We proceed by induction on $d$. Since the residue class field $A/\mathfrak{m}$ of $A$ is infinite, we may choose $a \in Q \backslash \mathfrak{m}Q$ is a superficial element for $I$. We set $\overline{A} = A/(a)$, $\overline{I} = /(a)$, $\overline{Q} = /(a)$, and $r' = \mathrm{r}_{\overline{Q}}(\overline{I})$.

Suppose that $d = 2$. We set $v_k = \ell_A(\overline{I}^{k+2}/\overline{Q}\,\overline{I}^{k+1})$ for all $k \geq 0$. Then, we have $g = \sum_{k=0}^{r'-2} v_k$ by Lemma 2.3. If $v_0 = \ell_A(\overline{I}^2/\overline{Q}\,\overline{I}) \leq 1$ then we have $\ell_A(I^2/QI) \leq 1$, because $\widetilde{I} = I$ by Lemma 4.3. Then we have $\mathrm{e}_2(I) = \binom{g+1}{2}$ by Theorem 1.2, but it is impossible. Therefore $v_0 = \ell_A(\overline{I}^2/\overline{Q}\,\overline{I}) \geq 2$. Then because $\sum_{k=1}^{r'-2} v_k = g - v_0 \leq g - 2$, we have

$$\sum_{k \geq 1} k{\cdot}v_k \leq \binom{g - v_0 + 1}{2} \leq \binom{g-1}{2}$$

by Lemma 3.1. Then, by Lemma 2.3, we have

$$\begin{aligned} \mathrm{e}_2(I) &= \sum_{k=0}^{r'-2}(k+1){\cdot}v_k - \sum_{k=1}^{m-1} \ell_A([I^{k+1} :_A a]/I^k) \\ &\leq \sum_{k=1}^{r'-2} k{\cdot}v_k + \sum_{k=0}^{r'-2} v_k \leq \binom{g-1}{2} + g = \binom{g}{2} + 1 \leq \mathrm{e}_2(I), \end{aligned}$$

where $m \geq 1$ be an integer such that $I^{k+1} :_A a = I^k$ for all $k \geq m$. Therefore, we have

$$\sum_{k \geq 1} k{\cdot}v_k = \binom{g - v_0 + 1}{2} = \binom{g-1}{2}$$

so that $v_0 = \ell_A(\overline{I}^2/\overline{Q}\,\overline{I}) = 2$, $v_k = \ell_A(\overline{I}^{k+2}/\overline{Q}\,\overline{I}^{k+1}) = 1$ for all $1 \leq k \leq g-2$, and $v_{g-1} = \ell_A(\overline{I}^{g+1}/\overline{Q}\,\overline{I}^g) = 0$ by Lemma 3.1. We also have $I^{k+1} :_A a = I^k$ for all $k \geq 1$, so that $\operatorname{depth} \operatorname{gr}(I) \geq 1$. Thus, we get $\ell_A(I^2/QI) = 2$, $\ell_A(I^{k+1}/QI^k) = 1$ for all $1 \leq k \leq g-1$, and $I^{g+1} = QI^g$ since $at$ forms a $\operatorname{gr}(I)$-regular element.

Assume that $d \geq 3$ and that our assertion holds true for $d-1$. Since $\mathrm{g}_s(\overline{I}) = g \geq 2$ by Lemma 2.2, we have

$$\binom{\mathrm{g}_s(\overline{I})}{2} + 1 \leq \mathrm{e}_2(\overline{I}) < \binom{\mathrm{g}_s(\overline{I})+1}{2}.$$

Thus, by the hypothesis of induction on $d$, we get $\ell_A(\overline{I}^2/\overline{Q}\,\overline{I}) = 2$, $\ell_A(\overline{I}^{k+1}/\overline{Q}\,\overline{I}^k) = 1$ for all $2 \leq k \leq g-1$, $\overline{I}^{g+1} = \overline{Q}\,\overline{I}^g$, and $\operatorname{depth} \operatorname{gr}(\overline{I}) \geq d-2 > 0$. Then, thanks to Sally's technique ([S1], [HM, Lemma 2.2]), the element $at$ is $\operatorname{gr}(I)$-regular. Therefore, we get $\ell_A(I^2/QI) = 2$, $\ell_A(I^{k+1}/QI^k) = 1$ for all $2 \leq k \leq g-1$, $I^{g+1} = QI^g$, and $\operatorname{depth} \operatorname{gr}(I) \geq d-1$ as required.

(3) $\Rightarrow$ (1) and $(i)$ We have $\sum_{k \geq 1} \ell_A(I^{k+1}/QI^k) = g$ by our assumption, so that the equality $\mathrm{e}_1(I) = \mathrm{e}_0(I) - \ell_A(A/I) + g = \sum_{k \geq 0} \ell_A(I^{k+1}/QI^k)$ holds true. Therefore, thanks to Theorem 2.4, we get $\mathrm{e}_2(I) = \binom{g}{2} + 1$ and $\mathrm{e}_i(I) = \binom{g}{i}$ for all $3 \leq i \leq d$.

$(iii)$ Assume that $\operatorname{gr}(I)$ is a Cohen-Macaulay ring. Then we have $Q \cap I^{n+1} = QI^n$ for all $n \geq 1$ ([VV]). Since $\ell_A(I^3/QI^2) \leq 1$, we have $I^4 \subseteq QI^2$, so that we get $I^4 = Q \cap I^4 = QI^3$. The converse also holds true by Valabrega-Valla criterion ([VV]). This completes the proof of Theorem 4.1. $\square$

In the end of this section we introduce one consequence of Theorem 4.1. The following result is the main theorem in [S2], which is exactly the case where $\mathrm{g}_s(I) = \mathrm{e}_2(I) = 2$ in Theorem 4.1.

**Corollary 4.4.** ([S2, Theorem 4.4]) *Suppose that $d \geq 2$. Assume that $\mathrm{g}_s(I) = \mathrm{e}_2(I) = 2$. Then the following assertions hold true.*

(1) $\mathrm{e}_i(I) = 0$ *for all* $3 \leq i \leq d$,
(2) $\ell_A(I^2/QI) = 2$ *and* $I^3 = QI^2$, *and*
(3) $\operatorname{depth} \operatorname{gr}(I) \geq d-1$.

## 5. Examples

In this section, we shall construct several examples of $\mathfrak{m}$-primary ideals $I$ satisfying the conditions of our main theorems (Theorem 3.4 and Theorem 4.1).

**Remark 5.1.** To construct necessary examples we may assume that $\dim A = 1$. In fact, let $d \geq 2$ be an integer and let $J$ be an $\mathfrak{n}$-primary ideal in a certain 1-dimensional Cohen-Macaulay local ring $(B, \mathfrak{n})$ with $\mathfrak{q} = aB$ a reduction of $J$. Let $A = B[[X_1, X_2, \cdots, X_{d-1}]]$ be the formal power series ring. We set $I = JA + (X_1, X_2, \cdots, X_{d-1})A$ and $Q = \mathfrak{q}A + (X_1, X_2, \cdots, X_{d-1})A$. Then $A$ is a Cohen-Macaulay local ring with $\dim A = d$ and the maximal ideal $\mathfrak{m} = \mathfrak{n}A + (X_1, X_2, \cdots, X_{d-1})A$. The ideal $Q$ is a reduction of $I$ and because $X_1, X_2, \cdots, X_{d-1}$ forms a super regular sequence for $I$ (recall that $\operatorname{gr}(I) = \operatorname{gr}(J)[Y_1, Y_2, \cdots, Y_{d-1}]$ is the polynomial ring, where $Y_i$'s are the initial forms of $X_i$'s), we have $\operatorname{depth} \operatorname{gr}(I) = \operatorname{depth} \operatorname{gr}(J) + d - 1 \geq d - 1$, $\mathrm{e}_i(I) = \mathrm{e}_i(J)$ $(i = 0, 1, 2)$, $\mathrm{g}_s(I) = \mathrm{g}_s(J)$, and $I^{k+1}/QI^k \cong J^{k+1}/\mathfrak{q}J^k$ for all $k \geq 0$. This observation allows us to concentrate our attention on the case where $\dim A = 1$.

We give the following example of $\mathfrak{m}$-primary ideals $I$ satisfying the equality $\mathrm{e}_2(I) = \binom{\mathrm{g}_s(I)+1}{2}$ in Theorem 3.4 (Theorem 1.1).

**Example 5.2.** Let $e \geq 3$ be an integer and let $\mathrm{H} = \langle 2e, 2e+1, \cdots, 4e-1 \rangle$ be the numerical semi-group generated by $2e, 2e+1, \cdots, 4e-1$. Let $B = k[[u^{2e}, u^{2e+1}, \cdots, u^{4e-1}]] \subseteq k[[u]]$, where $k[[u]]$ denotes the formal power series ring with one indeterminate $u$ over an infinite field $k$. Then $B$ is a 1-dimensional Cohen-Macaulay local ring with the maximal ideal $\mathfrak{n} = (u^{2e}, u^{2e+1}, \cdots, u^{4e-1})$. Let $\mathfrak{q} = (u^{2e})$ and $J = (u^{2e}, u^{2e+2}, u^{4e+1})$. Then we have the following.

- $(i)$ $\mathrm{g}_s(J) = e - 2$,
- $(ii)$ $\ell_B(J^{k+1}/\mathfrak{q}I^k) = 1$ for all $1 \leq k \leq e - 2$, and $I^e = QI^{e-1}$, and
- $(iii)$ $\operatorname{depth} \operatorname{gr}(J) = 0$.

Let $d \geq 2$ be an integer and $A = B[[X_1, X_2, \cdots, X_{d-1}]]$ be a formal power series ring over $B$, and $I = JA + (X_1, X_2, \cdots, X_{d-1})A$ and $Q = \mathfrak{q} + (X_1, X_2, \cdots, X_{d-1})A$. Then $Q$ is a reduction of $I$ and the following assertions hold true.

1. $A$ is a Cohen-Macaulay ring with $\dim A = d \geq 2$,
2. $\mathrm{g}_s(I) = e - 2$,
3. $\mathrm{e}_0(I) = \mathrm{e}_1(I) = e$ and $\mathrm{e}_i(I) = \binom{e-1}{i}$ for all $2 \leq i \leq d$, and
4. $\operatorname{depth} \operatorname{gr}(I) = d - 1$.

*Proof.* Thanks to Theorem 3.4 and Remark 5.1, we have only to show that the assertions $(i)$, $(ii)$, and $(iii)$ hold true.

It is routine to show that $J^{k+1} = \mathfrak{q}^k J + \sum_{\ell=1}^{k} \mathfrak{q}^{k-\ell} u^{2(\ell+1)(e+1)} = \mathfrak{q}J^k + u^{2(k+1)(e+1)}B$ for all $k \geq 1$. Then, since $u^{2k(e+1)} \notin \mathfrak{q}J^k$ for all $1 \leq k \leq e - 2$ and $\mathfrak{m}J \subseteq \mathfrak{q}$, we get $\ell_B(J^{k+1}/\mathfrak{q}J^k) = 1$ for all $1 \leq k \leq e - 2$ and $J^e = \mathfrak{q}J^{e-1}$. Therefore $\mathrm{g}_s(J) = \ell_B(B/J) - \mathrm{e}_0(I) + \mathrm{e}_1(I) = \sum_{k \geq 1} \ell_B(J^{k+1}/\mathfrak{q}J^k) = e - 2$ by Lemma 2.1. Since $J^2 \subseteq \mathfrak{q}$ and $J^2 \neq \mathfrak{q}J$, the associated graded ring $\operatorname{gr}(J)$ is not Cohen-Macaulay. $\square$

In the end of this paper, we introduce the example of $\mathfrak{m}$-primary ideals $I$ satisfying the equality $\mathrm{e}_2(I) = \binom{\mathrm{g}_s(I)}{2} + 1$ in Theorem 4.1 (Theorem 1.2).

**Example 5.3.** Let $e \geq 5$ be an integer and let $\mathrm{H} = \langle e, e+1, \cdots, 2e-1 \rangle$ be the numerical semi-group generated by $e, e+1, \cdots, 2e-1$. Let $d \geq 2$ be an integer and let $A = k[[u^e, u^{e+1}, \cdots, u^{2e-1}]] \subseteq k[[u]]$, where $k[[u]]$ denotes the formal power series ring with one indeterminate $u$ over an infinite field $k$. Then $A$ is a 1-dimensional Cohen-Macaulay local ring with the maximal ideal $\mathfrak{m} = (u^e, u^{e+1}, \cdots, u^{2e-1})$. Let $Q = (u^e)$ and $I = (u^e, u^{e+1}, u^{2e-2})$. We then have the following.

- $(i)$ $\mathrm{g}_s(J) = e - 3$,
- $(ii)$ $\ell_B(J^2/\mathfrak{q}J) = 2$, $\ell_B(J^{k+1}/\mathfrak{q}I^k) = 1$ for all $2 \leq k \leq e - 4$, and $J^{e-2} = \mathfrak{q}J^{e-3}$, and

(*iii*) $\operatorname{depth} \operatorname{gr}(J) = 0$.

Let $d \geq 2$ be an integer and $A = B[[X_1, X_2, \cdots, X_{d-1}]]$ be a formal power series ring over $B$, and $I = JA + (X_1, X_2, \cdots, X_{d-1})A$ and $Q = \mathfrak{q}A + (X_1, X_2, \cdots, X_{d-1})A$. Then $Q$ is a reduction of $I$ and the following assertions hold true.

(1) $A$ is a Cohen-Macaulay local ring with $\dim A = d \geq 2$,

(2) $\mathrm{g}_s(I) = e - 3$,

(3) $\mathrm{e}_0(I) = e$, $\mathrm{e}_1(I) = e - 1$, and $\mathrm{e}_i(I) = \binom{e-3}{i}$ for all $2 \leq i \leq d$, and

(4) $\operatorname{depth} \operatorname{gr}(I) = d - 1$.

*Proof.* It is routine to show that $J^2 = \mathfrak{q}J + (u^{2e+2}, u^{3e-1})B$ and $J^{k+1} = \mathfrak{q}^k J + \mathfrak{q}^{k-1}(u^{2e+2}, u^{3e-1}) + \sum_{\ell=2}^{k} \mathfrak{q}^{k-\ell} u^{(\ell+1)(e+1)} = \mathfrak{q}J^k + u^{(k+1)(e+1)}B$ for all $k \geq 2$. Therefore, since $u^{(k+1)(e+1)} \notin \mathfrak{q}J^k$ for all $2 \leq k \leq e-4$ and $\mathfrak{m}J \subseteq \mathfrak{q}$, we get $\ell_B(J^{k+1}/\mathfrak{q}J^k) = 1$ for all $2 \leq k \leq e-4$. We also have $\ell_B(J^2/\mathfrak{q}J) = 2$ and $J^{e-2} = \mathfrak{q}J^{e-3}$. Hence $\mathrm{g}_s(J) = \ell_B(B/J) - \mathrm{e}_0(I) + \mathrm{e}_1(I) = \sum_{k \geq 1} \ell_B(J^{k+1}/\mathfrak{q}J^k) = e - 3$ by Lemma 2.1. Since $J^2 \subseteq \mathfrak{q}$ and $J^2 \neq \mathfrak{q}J$, the associated graded ring $\operatorname{gr}(J)$ is not Cohen-Macaulay. Thus, thanks to Remark 5.1 and Theorem 4.1, all assertions (1), (2), (3), and (4) hold true. □

## References

[CPP] A. Corso, C. Polini, and M. V. Pinto, *Sally Modules and Associated Graded Rings*, Comm. Alg. **26** (1998), 2689–2708.

[HM] S. Huckaba and T. Marley, *Hilbert coefficients and the depth of associated graded rings*, J. London Math. Soc. (2) **56**, 1997, 64–76.

[H] C. Huneke, *Hilbert functions and symbolic powers*, Michigan Math. J., Vol **34**, 1987, 293–318.

[KM] D. Kirby and H. A. Mehran, *A note on the coefficients of the Hilbert-Samuel polynomial for a Cohen-Macaulay module*, J. London Math. Soc. (2) **25**, 1982, 449–457.

[Na] M. Narita, *A note on the coefficients of Hilbert characteristic functions in semi-regular rings*, Proc. Cambridge Philos. Soc. **59** (1963), 269–275.

[No] D. G. Northcott, *A note on the coefficients of the abstract Hilbert function*, J. London Math. Soc., Vol **35**, 1960, 209–214.

[O] A. Ooishi, *$\delta$-genera and sectional genera of commutative rings*, Hiroshima Math. J. **17**, 1987, 361–372.

[RR] L. J. Ratliff and D. Rush, *Two notes on reductions of ideals*, Indiana Univ. Math. J. **27**, 1978, 929-934.

[RV1] M. E. Rossi and G. Valla, *A conjecture of J. Sally*, Comm. Alg. **24** (1996), 4249–4261.

[RV2] M. E. Rossi and G. Valla, *Hilbert functions of filtered modules*, UMI Lecture Notes 9, Springer (2010).

[S1] J. D. Sally, *Hilbert coefficients and reduction number 2*, J. Alg. Geo. and Sing. **1**, 1992, 325–333.

[S2] J. D. Sally, *Ideals whose Hilbert function and Hilbert polynomial agree at $n = 1$*, J. Algebra **157** (1993), 534–547.

[VV] P. Valabrega and G. Valla, *Form rings and regular sequences*, Nagoya Math. J. **72** (1978), 93–101.

[W] H. J. Wang, *On Cohen-Macaulay local rings with embedding dimension $e + d - 2$*, J. Algebra **190** (1997), 226–240.

Department of Mathematical Science, Faculty of Science, Yamaguchi University, 1677-1 Yoshida, Yamaguchi 753-8512, Japan

*E-mail address*: ozeki@yamaguchi-u.ac.jp